\documentclass[12pt,a4paper]{article}
\usepackage[latin2]{inputenc}
\usepackage{amsmath}
\usepackage{amsfonts}
\usepackage{amssymb}
\usepackage{graphicx}
\usepackage[left=2.00cm, right=2.00cm, top=2.50cm, bottom=2.50cm]{geometry}
\usepackage[T1]{fontenc}

\def\Q{{\mathbb Q}}
\def\Z{{\mathbb Z}}

\newtheorem{lemma}{Lemma}
\newtheorem{theorem}[lemma]{Theorem}

\title{
Totally real Thue inequalities \\ over imaginary quadratic fields: an improvement
}
\author{
Istv\'{a}n Ga\'{a}l\thanks{
        Research supported in part by  the EFOP-3.6.1-16-2016-00022 project. 
				The project is co-financed by the European Union and the European Social Fund.	
                         },\; \\
{\small University of Debrecen, Mathematical Institute} \\
{\small H--4002 Debrecen Pf.400., Hungary,} 
{\small e--mail: gaal.istvan@unideb.hu},
 \\ \\
Borka Jadrijevi\'c\\
{\small University of Split,}
{\small Faculty of Science, }\\
{\small Ru\dj era Bo\v{s}kovi\'{c}a 33, 21000 Split, Croatia,}
{\small e--mail: borka@pmfst.hr}
 \\ \\ 
L\'aszl\'o Remete\thanks{
        Research supported by the \'UNKP-19-3 new national excellence program of the Ministry of human capacities.}\; \\
{\small University of Debrecen, Mathematical Institute} \\
{\small H--4002 Debrecen Pf.400., Hungary,} 
{\small e--mail: remete.laszlo@science.unideb.hu}
}

\begin{document}
\baselineskip=17pt

\maketitle
\thispagestyle{empty}

\renewcommand{\thefootnote}{\arabic{footnote}}
\setcounter{footnote}{0}


\begin{abstract}
We significantly improve our results of \cite{gjr} reducing relative 
Thue inequalities to absolute ones.
\end{abstract}

\section{Results}

Let $F(x,y)$ be a binary form of degree $n\geq 3$ with rational integer coefficients.
Assume that $f(x)=F(x,1)$ has leading coefficient 1 and 
distinct real roots $\alpha_1,\ldots,\alpha_n$. 
Let $0<\varepsilon<1$ and let $K\geq 1$. 
Let
\[
A=\min_{i\not=j}|\alpha_i-\alpha_j|,\;\; 
B=\min_i\prod_{j\not=i} |\alpha_j-\alpha_i|,\;\;
C=\frac{K}{(1-\varepsilon)^{n-1}B},\;\;
G=\frac{K^{1/n}}{\varepsilon A}.
\]

\vspace{1cm}

\noindent
Let $m\ge 1$ be a square-free positive integer, and set $M=\Q(i\sqrt{m})$. 
Consider the relative inequality
\begin{equation}
|F(x,y)|\leq K \;\; {\rm in}\;\; x,y\in\Z_M.
\label{1}
\end{equation}
If $F$ is irreducible, then (\ref{1}) is called a Thue inequality. We emphasize that
our statements are valid also if $F$ is reducible.

\vspace{1cm}

\noindent
If $m\equiv 3\; (\bmod \; 4)$, then $x,y\in\Z_M$ can be written as
\[
x=x_1+x_2\frac{1+i\sqrt{m}}{2}=\frac{(2x_1+x_2)+x_2i\sqrt{m}}{2}, 
\;\;
y=y_1+y_2\frac{1+i\sqrt{m}}{2}=\frac{(2y_1+y_2)+y_2i\sqrt{m}}{2},
\]
and if $m\equiv 1,2\; (\bmod \; 4)$, then 
\[
x=x_1+x_2i\sqrt{m},\;\; y=y_1+y_2i\sqrt{m},
\]
in both cases with $x_1,x_2,y_1,y_2\in\Z$.
Set
$s=2$ if $m\equiv 3\; (\bmod\;4)$ and $s=1$ if $m\equiv 1,2\; (\bmod\;4)$.
In the following Theorem  we formulate our statements parallely in the two cases.

\begin{theorem}
\label{th1}
Let $(x,y)\in\Z_M^2$ be a solution of (\ref{1}).  Then\\
\begin{equation}
|F(sx_1+(s-1)x_2,sy_1+(s-1)y_2)|\leq s^n K,\;\; |F(x_2,y_2)|\leq \frac{s^nK}{(\sqrt{m})^n},
\label{a}
\end{equation}
\begin{equation}
{\it and}\;\;\;\; |F(sx_1+(s-1)x_2,sy_1+(s-1)y_2)|\cdot |F(x_2,y_2)|\leq \frac{s^{2n}K^2}{2^n\cdot (\sqrt{m})^n}.
\label{aa}
\end{equation}
\begin{equation}
{\rm If}\;\; |y|>\max\left\{G,
\left(\frac{s\cdot C}{\sqrt{m}}\right)^{\frac{1}{n-2}}\right\},
\;\; {\rm then} \;\; x_2y_1=x_1y_2.
\label{x12}
\end{equation}
\begin{equation}
{\rm If}\;\; |y|>\max\left\{G, (s\cdot C)^{\frac{1}{n-1}} \right\}
\;\;{\rm and}\;\; sy_1+(s-1)y_2=0, \;\;{\rm then}\;\; sx_1+(s-1)x_2=0.
\label{I1}
\end{equation}
\begin{equation}
{\rm If}\;\; |y|>\max\left\{G, \left(\frac{s\cdot C}{\sqrt{m}}\right)^{\frac{1}{n-1}} \right\}
\;\;{\rm and}\;\; y_2=0, \;\;{\rm then}\;\; x_2=0.
\label{I2}
\end{equation}
\end{theorem}

\vspace{0.5cm}

\noindent
{\bf Remark 1.} The present inequality (\ref{a})  is
much sharper than the corresponding inequalities of Theorem 2.1 of \cite{gjr}. Moreover
we obtain these inequalities without any conditions on the variables. This makes the 
applications much easier. If the values of $F$ are non-zero, then (\ref{aa}) 
yields further new restrictions for the possible solutions of (\ref{1}).

\vspace{0.5cm}

\noindent
{\bf Proof of Theorem \ref{th1}.}\\
Let $(x,y)\in\Z_M^2$ be an arbitrary solution of (\ref{1}).
Let $\beta_j=x-\alpha_j y,\; j=1,\ldots,n$, then inequality (\ref{1}) can be written as
\begin{equation}
|\beta_1\cdots\beta_n|\leq K.
\label{e}
\end{equation}
We have
\[
\beta_{j}=\frac{1}{s}((sx_{1}+\left(  s-1\right)  x_{2})-\alpha_{j}%
(sy_{1}+\left(  s-1\right)  y_{2}))+\frac{i\sqrt{m}}{s}(x_{2}-\alpha_{j}%
y_{2}).
\]
Obviously,
\[
|{\rm Re}(\beta_j)|\leq |\beta_j|,\;\; |{\rm Im}(\beta_j)|\leq |\beta_j|,\; 1\leq j\leq n.
\]
Further,
\[
\prod_{j=1}^n|{\rm Re}(\beta_j)|\leq \prod_{j=1}^n|\beta_j|\leq K,
\;\;
{\rm and}
\;\;
\prod_{j=1}^n|{\rm Im}(\beta_j)|\leq \prod_{j=1}^n|\beta_j|\leq K,
\]
which imply (\ref{a}). Moreover,
\[
\prod_{j=1}^n|{\rm Re}(\beta_j)|\cdot \prod_{j=1}^n|{\rm Im}(\beta_j)|=
\prod_{j=1}^n\left(|{\rm Re}(\beta_j)|\cdot |{\rm Im}(\beta_j)|\right)
\le
\prod_{j=1}^n \frac{|{\rm Re}(\beta_j)|^2+|{\rm Im}(\beta_j)|^2}{2}
= \prod_{j=1}^n \frac{|\beta_j|^2}{2}
\le \frac{K^2}{2^n},
\]
whence we obtain (\ref{aa}).

Assume now 
\begin{equation}
|y|\geq G.
\label{bbqq}
\end{equation}
Let ${i_0}$ be the index with
$|\beta_{i_0}|=\min_j |\beta_j|$.
Then $|\beta_{i_0}|\leq K^{\frac{1}{n}}$ and for $j\not={i_0}$
\begin{equation}
|\beta_j|\geq |\beta_j-\beta_{i_0}|-|\beta_{i_0}|\geq |\alpha_j-\alpha_{i_0}|\cdot |y|-K^{\frac{1}{n}}
\geq (1-\varepsilon) \cdot |\alpha_j-\alpha_{i_0}|\cdot |y|.
\label{bj}
\end{equation}
From (\ref{e}) and (\ref{bj}) we have
\begin{equation}
|\beta_{i_0}|\leq \frac{K}{\prod_{j\not={i_0}}|\beta_j|}\leq \frac{C}{|y|^{n-1}}.
\label{betai}
\end{equation}
Using that $\alpha_{i_0}|y|^2$ is real, by (\ref{betai}) we obtain
\[
|{\rm Im} (x\overline{y})|=|{\rm Im} (\alpha_{i_0}|y|^2-x\overline{y})|
\leq 
\left|\alpha_{i_0}|y|^2-x\overline{y}\right|
=
|y|^2\cdot \left|\alpha_{i_0}-\frac{x\overline{y}}{y\overline{y}}\right|
=
|y|^2\cdot \left|\alpha_{i_0}-\frac{x}{y}\right|
\leq \frac{C}{|y|^{n-2}}.
\]
If 
\[
|y|>\left(\frac{s\cdot C}{\sqrt{m}}\right)^{\frac{1}{n-2}},
\]
then this implies $x_2y_1=x_1y_2$.

Inequality (\ref{betai}) indicates that $|\beta_{i_0}|$ is small for sufficiently
large $|y|$ and so are its real and imaginary parts that can equal zero if we impose some
extra assumptions.

--If $|y|>(sC)^{\frac{1}{n-1}}$, then $|(sx_{1}+\left(  s-1\right)
x_{2})-\alpha_{i_{0}}(sy_{1}+\left(  s-1\right)  y_{2})|<1$. So,
$sy_{1}+\left(  s-1\right)  y_{2}=0$ implies $sx_{1}+\left(  s-1\right)
x_{2}=0$.

--If $|y|>\displaystyle{\left(  \frac{sC}{\sqrt{m}}\right)  ^{\frac{1}{n-1}}}%
$, then $|x_{2}-\alpha_{i_{0}}y_{2}|<1$. So, $y_{2}=0$ implies $x_{2}=0$.
\hfill $\Box$

\section{How to apply Theorem \ref{th1}?}

Finally, we give useful hints for a practical application of Theorem \ref{th1}.
Using the same notation let us consider again the relative 
inequality (\ref{1}). We describe our algorithm in case $m\equiv 3\;(\bmod\;4)$,
since the other case is completely similar.\\

\noindent
First, we solve $F(x_{2},y_{2})=k_1$ for all $k_1\in\Z$ with
$|k_1|\leq 2^{n}K/ (\sqrt{m})^n.$ 
Since the equation $F(x_{2},y_{2})=0$ can also have non-trivial solutions
if $F$ is reducible, we split our arguments into two cases.

A. First suppose $F(x_{2},y_{2})=0$. This makes possible to determine $x_2,y_2$.
If $F$ is irreducible, then $x_2=y_2=0$, if $F$ is reducible, then $x_2,y_2$ can be
determined easily (if there are any). 
We then determine the solutions $(a,b)\in\Z^2$ of $|F(a,b)|=k_2$ 
for all $k_2$ with $|k_2|\le 2^n K$.
Using all possible values of $x_2,y_2$ for each solution $(a,b)$ we determine
$x_1=(a-x_2)/2, y_1=(b-y_2)/2$ and check if these are integers.
(Note that if $F$ is irreducible, then $x_2=y_2=0$ implies $|F(x_1,y_1)|\le K$
and the procedure can be simplified.)
Having all possible $x_1,x_2,y_1,y_2$ we test if 
$(x,y)\in\Z_M^2$ is a solution of (\ref{1}).

B. Assume now $F(x_{2},y_{2})=k_1\ne 0$ for some $(x_2,y_2)\in\Z^2$. 
Then we solve $F(a,b)=k_2$ in $(a,b)\in\Z^2$ for all
$k_2\in\Z$ with $|k_1k_2|\le 2^{n}K^2/ (\sqrt{m})^n$
(a part of this calculation was already performed by solving 
$F(x_{2},y_{2})=k_1$). Having $a,b,x_2,y_2$ we calculate 
$x_1=(a-x_2)/2,y_1=(b-y_2)/2$. For $x_2,y_2$ and integer values $x_1,y_1$ we test if 
$(x,y)^2\in\Z_M$ is indeed a solution of (\ref{1}).

\vspace{0.5cm}

\noindent
{\bf Remark 2.}
If $m$ is sufficiently large, then by (\ref{a}) we have $|F(x_2,y_2)|<1$.
In case $F$ is irreducible, this implies $x_2=y_2=0$, whence (\ref{1}) reduces to
an inequality in $x_1,y_1$ over $\Z$.

\vspace{0.5cm}

\noindent
{\bf Remark 3.}
Solving Thue equations over $\Z$ is no problem any more by using
well-known computer algebra packages. If $F$ is reducible, this task is even 
easier.

\end{document}